\documentclass[12pt,fullpage]{article}

\usepackage{amsfonts}
\usepackage{xr}
\usepackage{graphicx}
\usepackage{amsmath,amssymb}
\usepackage{dsfont,epsfig,subfigure}
\usepackage{color}
\usepackage{epstopdf}
\def\R{\mathbb{R}}
\def\N{\mathbb{N}}

\def\epsilon{\varepsilon}

\def\tilde{\widetilde}

\def\essinf{{\rm{ess}}\ {\rm{inf}}}
\def\esssup{{\rm{ess}}\ {\rm{sup}}}

\newcommand{\SE}{\setcounter{equation}{0} \section}
\newcommand{\be}{\begin{equation}}
\newcommand{\ee}{\end{equation}}
\newcommand{\baa}{\begin{array}}
\newcommand{\eaa}{\end{array}}
\newcommand{\ba}{\begin{eqnarray}}
\newcommand{\ea}{\end{eqnarray}}

\newtheorem{theo}{\bf Theorem}[section]

\newtheorem{cor}[theo]{\bf Corollary}
\newtheorem{defi}[theo]{\bf Definition}
\newtheorem{rem}[theo]{\bf Remark}

\newtheorem{conj}[theo]{\bf Conjecture}

\begin{document}
\date{}
\title{\bf{Transition fronts and stretching phenomena for a general class of reaction-dispersion equations}}
\author{Jimmy Garnier$^{\hbox{\small{ a}}}$, Fran\c cois Hamel$^{\hbox{\small{ b}}}$ and Lionel Roques$^{\hbox{\small{ c}}}$\thanks{This work has been carried out in the framework of Archim\`ede LabEx (ANR-11-LABX-0033) and of the A*MIDEX project (ANR-11-IDEX-0001-02), funded by the ``Investissements d'Avenir" French Government program managed by the French National Research Agency (ANR). The research leading to these results has received funding from the European Research Council under the European Union's Seventh Framework Programme (FP/2007-2013) / ERC Grant Agreement n.321186 - ReaDi - Reaction-Diffusion Equations, Propagation and Modelling, and from the ANR project NONLOCAL (ANR-14-CE25-0013).}\\
\\
\footnotesize{$^{\hbox{a }}$Universit\'e de Savoie, Laboratoire de Math\'ematiques (LAMA) - EDPs$^2$ (UMR 5127)}\\
\footnotesize{B\^at. Le Chablais, Campus Scientifique, 73376 Le Bourget du Lac, France}\\
\footnotesize{$^{\hbox{b }}$Aix Marseille Universit\'e, CNRS, Centrale Marseille}\\
\footnotesize{Institut de Math\'ematiques de Marseille, UMR 7373, 13453 Marseille, France}\\
\footnotesize{$^{\hbox{c }}$INRA, UR546 Biostatistique et Processus Spatiaux, 84914 Avignon, France}}
\maketitle

\begin{center}
{\it{Dedicated to the memory of Professor Paul Fife}}
\end{center}
\vskip 0.3cm

\begin{abstract}
 We consider a general form of reaction-dispersion equations with non-local dispersal and local reaction. Under some general conditions, we prove the non-existence of transition fronts, as well as some stretching properties at large time for the solutions of the Cauchy problem. These conditions are satisfied in particular when the reaction is monostable and when the dispersal operator is either the fractional Laplacian, a convolution operator with a fat-tailed kernel or a nonlinear fast diffusion operator.
\end{abstract}


\SE{Introduction}\label{intro}

This note is concerned with { fast } propagation phenomena and large time qualitative pro\-perties of dispersion-reaction equations of the type
\be\label{eq}
u_t(t,x)=\mathcal{D}u(t,x)+f(u(t,x)),\ \ t\in\R,\ x\in\R.
\ee
The given reaction function $f:[0,1]\to\R$ is of class $C^1$ and such that $f(0)=f(1)=0$. Throughout the paper, we assume that, for any $u_0\in L^{\infty}(\R;[0,1])$, the Cauchy problem
\be\label{cauchy}\left\{\baa{rcll}
u_t(t,x) & = & \mathcal{D}u(t,x)+f(u(t,x)), & t>0,\,x\in\R,\vspace{3pt}\\
u(0,\cdot) & = & u_0 & \eaa\right.
\ee
admits a unique mild solution, such that $u(t,\cdot)\in L^{\infty}(\R;[0,1])$ for every $t>0$ and $u$ is uniformly continuous with respect to $t$ in $[\delta,+\infty)\times\R$ for every $\delta>0$. The equation~\eqref{eq} is assumed to be autonomous in the sense that, for any solution $u$ of the Cauchy problem~\eqref{cauchy} and for any $T>0$, the function $(t,x)\mapsto u(T+t,x)$ solves~\eqref{cauchy} with initial condition $u(T,\cdot)$. We also assume that a comparison principle holds for~\eqref{cauchy}, that is, for any two solutions $u$ and $v$ with respective initial conditions $u_0$ and $v_0$,
\be\label{comp}
\big(0\le u_0\le v_0\le 1\hbox{ a.e. in }\R\big)\ \Longrightarrow\ \big(\forall\,t>0,\ 0\le u(t,\cdot)\le v(t,\cdot)\le 1\hbox{ a.e. in }\R\big).
\ee
Lastly, the solutions $u$ of~\eqref{cauchy} are assumed to have infinite spreading speeds (to the right) in the sense that
\be\label{spreading}
\big(u_0\in L^{\infty}(\R;[0,1])\backslash\{0\}\big)\ \Longrightarrow\ \Big(\forall\,c>0,\ \mathop{{\rm{ess}}\,{\rm{inf}}}_{(0,ct)}u(t,\cdot)\to1\hbox{ as }t\to+\infty\Big).
\ee

The first typical example of equations~\eqref{eq} for which these conditions are fulfilled is { the monostable reaction-diffusion equation with fractional diffusion }
\be\label{fraclaplace}
u_t=-(-\Delta_x)^{\alpha}u+f(u),
\ee
where the reaction function $f$ { is monostable, that is}
\be\label{f}
f(0)=f(1)=0,\ f(u)>0\hbox{ for all }u\in(0,1),\ f'(0)>0,
\ee
and the dispersal operator $\mathcal{D}$ is the fractional Laplacian $\mathcal{D}u=-(-\Delta_x)^{\alpha}u$ with $0\!<\!\alpha\!<\!1$, see \cite{cr} (see also~\cite{ccr} for similar equations in periodic media). In~\cite{cr}, the function~$f$ is also assumed to be concave in $[0,1]$, but the condition~\eqref{spreading} actually holds under the condition~\eqref{f} by the comparison principle and putting below $f$ a concave function $g:[0,1-\epsilon]\to\R$ on the interval~$[0,1-\epsilon]$, with $\epsilon\in(0,1)$ as small as wanted. The set of functions~$f$ satis\-fying~\eqref{f} includes the class of Fisher-KPP~\cite{fi,kpp} nonlinearities~$f$ for which, in addition to these assumptions,~$f$ satisfies $f(u)\le f'(0)u$ in $[0,1]$. It also contains the set of $C^2$ concave functions~$f$ with $f(0)=f(1)=0$ and~$f>0$ on $(0,1)$, whose archetype is the logistic nonlinearity~$f(u)=u(1-u)$.

A second important class consists of { integro-differential equations with fat-tailed dispersal kernels}
\be\label{convolution}
u_t=J*u-u+f(u),
\ee
where $f$ is of the type~\eqref{f} and the dispersal operator is a convolution operator
$$\mathcal{D}u(t,x)=J*u(t,x)-u(t,x):=\int_{\R}J(x-y)(u(t,y)-u(t,x)) \, dy,$$
with a fat-tailed positive even kernel $J\in L^1(\R)$ such that
$$\left\{\baa{l}
\displaystyle\|J\|_{L^1(\R)}=1,\ \int_{\R}|x|J(x)dx<+\infty,\vspace{3pt}\\
J\hbox{ is of class }C^1\hbox{ for large }|x|\hbox{ and }J'(x)=o(J(x))\hbox{ as }|x|\to+\infty,\eaa\right.$$
see~\cite{g}. Notice that such kernels~$J$ are called ``fat tailed" since they decay slowly as $|x|\to+\infty$ in the sense that $J(x)e^{\eta|x|}\to+\infty$ as $|x|\to+\infty$ for every $\eta>0$. Archetypes are $J(x)=c_{\alpha,\beta}e^{-\beta|x|^{\alpha}}$ with $0<\alpha<1$, $\beta>0$ and some normalization constant~$c_{\alpha,\beta}>0$, or $J(x)=c_{\alpha}/(1+|x|^{\alpha})$ with $\alpha>2$ and some $c_{\alpha}>0$. Such operators and equations arise in many physical or ecological models, see e.g.~\cite{cl,dcn,dcncl,f,klv,mvv,mk,mo}.

Lastly, the conditions \eqref{comp} and~\eqref{spreading} are also fulfilled when the dispersal operator $\mathcal{D}$ corresponds to fast nonlinear diffusion. Using a detailed formal analysis of the equation
\be
\label{eq:fastdif}
u_t=(u^\gamma)_{xx}+f(u),
\ee
with $0<\gamma<1$ and where $f$ satisfies \eqref{f}, it was shown in \cite{KinMcc03} that the assumption ~\eqref{spreading} was fulfilled. In \cite{StaVaz14}, this result was extended
to the case of nonlinear fractional diffusion equations:
\be
\label{eq:fastdif_frac}
u_t=-(-\Delta_x)^{\alpha}(u^\gamma)+f(u),
\ee
with $0<\alpha<1$ and $\max(1-2\, \alpha,0)<\gamma$ and where $f$ satisfies \eqref{f} and is a concave function. Due to the comparison principle~\eqref{comp}, which is valid for both~\eqref{eq:fastdif} and~\eqref{eq:fastdif_frac}, it can be shown that the condition~\eqref{spreading} holds without the concavity assumption, as for the standard fractional Laplacian.

Given the few general assumptions~\eqref{comp} and~\eqref{spreading}, the goal of the paper is twofold. Firstly, we will prove the non-existence of front-like entire solutions. Secondly, further stretching properties of the solutions of the Cauchy problem~\eqref{cauchy} at large time will be shown. As a corollary of the second main theorem, we will show that the solutions of~\eqref{cauchy} are in some sense flat at some large times in left and right neighborhoods of any level set. We again insist on the fact that these results will hold for the four main examples~\eqref{fraclaplace},~\eqref{convolution},~\eqref{eq:fastdif} and~\eqref{eq:fastdif_frac}.


\SE{Non-existence of transition fronts}

When propagating solutions are mentioned, one immediately has in mind standard trave\-ling fronts $u(t,x)=\varphi(\pm x-ct)$, with velocity $c\in\R$ and front profile $\varphi:\R\to[0,1]$ such that~$\varphi(-\infty)=1$ and $\varphi(+\infty)=0$. On the one hand, without the assumption~\eqref{spreading}, standard traveling fronts~$\varphi_c(\pm x-ct)$ are known to exist when we consider equations~\eqref{eq:fastdif} with the standard (local) Laplacian $\mathcal{D}u=\Delta_x u=u_{xx}$ or a nonlocal convolution operator~$\mathcal{D}u=J*u-u$ with thin-tailed dispersal kernel $J\in L^1(\R),$ which is nonnegative even and exponentially bounded in the sense that $\|J\|_{L^1(\R)}=1$ and $\int_{\R}J(x)\,e^{\lambda x}dx<+\infty$ for some $\lambda>0$. More precisely, when~$f$ is of the type~\eqref{f}, these fronts exist for every speed $c\ge c^*$, with a positive minimal speed~$c^*$. Furthermore, uniqueness of the profile~$\varphi_c$ and stability results for a given $c\ge c^*$ have been shown, with possibly heterogeneous nonlinearities~$f(x,u)$, see~\cite{aw,cc,co,cdm1,cd1,cd2,kpp,sa,s,sz1,sz2,u}. However, even for the reaction-diffusion equation $u_t=u_{xx}+f(u)$ with local diffusion and~$f$ of the type~\eqref{f}, fast propagation phenomena with infinite spreading speed as in~\eqref{spreading} is known to occur when $u_0$ decays to $0$ as $x\to\pm\infty$ more slowly than any exponentially decaying function, see~\cite{hr1}. On the other hand, some existence, uniqueness and stability results of standard traveling fronts have also been shown for some nonlocal equations of the type~\eqref{eq} with fractional Laplacians or fat-tailed dispersal convolution kernels $J$, but for other types of nonlinearities $f(u)$ or $f(x,u)$ than~\eqref{f}, see e.g.~\cite{bfrw,c,ch,co,cdm2,cd1,cd2,gh,gz,mrs2}.

Here, the assumptions~\eqref{comp} and~\eqref{spreading} forbid standard traveling fronts to exist, since all non-trivial solutions of the Cauchy problem accelerate with infinite speed as $t\to+\infty$ in the sense of~\eqref{spreading}. This non-existence result has already been known in the case of the fractional Laplacian $\mathcal{D}u=-(-\Delta_x)^{\alpha}u$~\cite{cr,gh} as well as in the convolution case $\mathcal{D}u=J*u-u$ with non-exponentially-bounded kernels $J$~\cite{g,y}. Non-existence of standard traveling fronts has also been derived for the fast nonlinear diffusion equations \eqref{eq:fastdif} and \eqref{eq:fastdif_frac} \cite{KinMcc03,StaVaz14}. However, other propagating solutions connecting $0$ and $1$, more general than the standard traveling fronts, can be investigated. These solutions, called transition fronts, have been defined in~\cite{bh4,bh5} in more general situations. In the one-dimensional situation considered here, as a wave-like solution defined in~\cite{s4}, the following definition holds.

\begin{defi}[Transition front]
A transition front connecting~$0$ $($say, on the right$)$ and~$1$ $($on the left$)$ for~\eqref{eq} is a time-global solution~$u:\R\times\R\to[0,1]$ $($that is, the unique mild solution of the Cauchy problem~\eqref{cauchy} starting at any time $t_0\in\R$ with ``initial" condition $u(t_0,\cdot))$ such that there exists a family~$(\xi_t)_{t\in\R}$ of real numbers with
\be\label{tf}\left\{\baa{ll}
\displaystyle\mathop{\essinf}_{(-\infty,x)}\,u(t,\xi_t+\cdot)\to 1 & \hbox{as }x\to-\infty,\vspace{3pt}\\
\displaystyle\mathop{\esssup}_{(x,+\infty)}\,u(t,\xi_t+\cdot)\to 0 & \hbox{as }x\to+\infty,\eaa\right.\hbox{ uniformly in }t\in\R.
\ee
\end{defi}

The real numbers $\xi_t$ therefore reflect the positions of a transition front $u$ as a function of time. However, the $\xi_t$'s are not uniquely defined since any family $(\tilde{\xi}_t)_{t\in\R}$ satisfies~\eqref{tf} as soon as~$(\xi_t)_{t\in\R}$ does and $(\tilde{\xi}_t-\xi_t)_{t\in\R}$ remains bounded. Roughly speaking, condition~\eqref{tf} means that the diameter of the transition zone between the sets where $u\simeq 1$ and $u\simeq 0$ is uniformly bounded in time in the sense that, for every $0<a\le b<1$, there is a nonnegative real number~$C$ such that, for all $t\in\R$,
$$\big\{x\in\R;\ a\le u(t,x)\le b\big\}\subset\big[X(t)-C,X(t)+C]$$
up to a negligible set. Transition fronts connecting $0$ on the left and $1$ on the right can be defined similarly by permu\-ting the limits as~$x\to\pm\infty$ in~\eqref{tf}. We will therefore only consider here transition fronts connec\-ting~$0$ on the right and $1$ on the left. Obviously, any standard traveling front~$u(t,x)=\varphi(x-ct)$ with~$\varphi(-\infty)=1$ and $\varphi(+\infty)=0$ would be a transition front connecting~$0$ and~$1$ with~$\xi_t=ct$, but, as already emphasized, the assumption~\eqref{spreading} excludes the existence of standard trave\-ling fronts for~\eqref{eq}. Nevertheless, transition fronts can a priori be much more general than standard traveling fronts, since no assumption is made on the family of front positions $(\xi_t)_{t\in\R}$. In particular, transition fronts different from the standard traveling fronts have been constructed recently for some homogeneous or heterogeneous local reaction-diffusion equations~\cite{hr2,hr3,mnrr,mrs1,nr1,nrrz,nr2,z1,z2}, and fronts with global speed $\xi_t/t$ having different limits as~$t\to\pm\infty$ are also known to exist even for the homogeneous local Fisher-KPP reaction-diffusion equation $u_t=u_{xx}+f(u)$~\cite{hn,hr2,z1}. For our equation~\eqref{eq}, transition fronts with global speed~$\xi_t/t$ growing arbitrarily as~$t\to\pm\infty$ are a priori not excluded.

Our first result shows actually that this is anyway impossible, since transition fronts cannot exist whatever the family of positions $(\xi_t)_{t\in\R}$ may be.

\begin{theo}\label{th1}
Under the assumptions~\eqref{comp} and~\eqref{spreading}, equation~\eqref{eq} does not admit any transition front connecting~$0$ and $1$.
\end{theo}

Theorem~\ref{th1} generalizes the known non-existence result of standard traveling fronts for the equations of the type~\eqref{fraclaplace} with the fractional Laplacian,~\eqref{convolution} with a fat-tailed dispersal kernel $J$,and \eqref{eq:fastdif} and \eqref{eq:fastdif_frac} with fast nonlinear diffusion since the standard fronts are particular classes of transition fronts.

Roughly speaking, this non-existence result can be heuristically explained as follows: for any transition transition front connecting $0$ and $1$, the uniformity (with respect to time) of the limits~\eqref{tf} prevents the front from traveling too fast (see also~\cite{hr2} for a similar phenomenon for local reaction-diffusion equations), and this last property is finally in contradiction with the spreading properties~\eqref{spreading} combined with the comparison principle~\eqref{comp}.\hfill\break

\noindent{\bf{Proof of Theorem~\ref{th1}.}} We argue by contradiction. So assume that, under the assumptions~\eqref{comp} and~\eqref{spreading}, equation~\eqref{eq} admits a transition front $u$ connecting $0$ and $1$. Let~$(\xi_t)_{t\in\R}$ be a family of real numbers such that~\eqref{tf} holds.\par
For every $t\in\R$, denote
$$\xi_t^-=\sup\,\Big\{x\in\R;\ \mathop{\essinf}_{(-\infty,x)}u(t,\cdot)\ge\frac{2}{3}\Big\}$$
and
$$\xi_t^+=\inf\,\Big\{x\in\R;\ \mathop{\esssup}_{(x,+\infty)}u(t,\cdot)\le\frac{1}{3}\Big\}.$$
It follows from~\eqref{tf} that $\xi_t^{\pm}$ are real numbers. The definitions of $\xi_t^{\pm}$ imply that
\be\label{xi-t}
\mathop{\essinf}_{(-\infty,\xi_t^-)}u(t,\cdot)\ge\frac23,\ \ \mathop{\essinf}_{(\xi_t^-,\xi_t^-+1)}u(t,\cdot)<\frac23
\ee
and
\be\label{xi+t}
\mathop{\esssup}_{(\xi_t^+,+\infty)}u(t,\cdot)\le\frac13,\ \ \mathop{\esssup}_{(\xi_t^+-1,\xi_t^+)}u(t,\cdot)>\frac13
\ee
together with $\xi_t^-\le\xi_t^+$. Property~\eqref{tf} also yields
$$\sup_{t\in\R}\big(\xi_t-\xi_t^-\big)<+\infty.$$
Indeed, otherwise, there would exist a sequence $(t_n)_{n\in\N}$ of real numbers such that
$$\xi_{t_n}-\xi^-_{t_n}\to+\infty\hbox{ as }n\to+\infty,$$
whence $\essinf_{(-\infty,\xi^-_{t_n}+1)}u(t_n,\cdot)\to1$ as $n\to+\infty$, contradicting the property
$$\essinf_{(\xi^-_{t_n},\xi^-_{t_n}+1)}u(t_n,\cdot)<\frac23.$$
Similarly, there holds
$$\inf_{t\in\R}\big(\xi_t-\xi_t^+\big)>-\infty.$$
Eventually, using again $\xi_t^-\le\xi_t^+$, one infers that the families $(\xi_t-\xi_t^{\pm})_{t\in\R}$ and $(\xi_t^+-\xi_t^-)_{t\in\R}$ are all bounded.\par
From the general regularity assumptions made in the paper, the function $u$ is actually uniformly continuous with respect to $t$ in $[T,+\infty)\times\R$, for every $T\in\R$. By considering $T=0$, it follows in particular that there exists $\delta>0$ such that, for all $t\ge0$, $s\ge0$ and~$x\in\R$,
$$\big(|t-s|\le\delta\big)\ \Longrightarrow\ \Big(|u(t,x)-u(s,x)|\le\frac{1}{12}\Big).$$
As a consequence, for every $t\ge0$, denoting $A_t^{\pm}$ the non-negligible measurable sets defined as
$$A^-_t=\Big\{x\in(\xi^-_t,\xi^-_t+1);\ u(t,x)<\frac{2}{3}\Big\}\ \hbox{ and }\ A^+_t=\Big\{x\in(\xi^+_t-1,\xi^+_t);\ u(t,x)>\frac{1}{3}\Big\},$$
it follows that
\be\label{134}
u(t+\delta,x)<\frac34\hbox{ for all }x\in A^-_t\ \hbox{ and }\ u(t+\delta,x)>\frac14\hbox{ for all }x\in A^+_t.
\ee\par
On the other hand, from~\eqref{tf}, there is $M\ge0$ such that, for every $t\in\R$,
$$\left\{\baa{ll}
\displaystyle u(t,x)\ge\frac34 & \hbox{for a.e. }x\in(-\infty,\xi_t-M),\vspace{3pt}\\
\displaystyle u(t,x)\le\frac14 & \hbox{for a.e. }x\in(\xi_t+M,+\infty).\eaa\right.$$
Since, for every $t\ge0$, the sets $A_t^{\pm}$ have a positive measure, one infers from~\eqref{134} that
$$\xi_t^-+1>\xi_{t+\delta}-M$$
and
$$\xi_t^+-1<\xi_{t+\delta}+M.$$
Remembering that the quantities $(\xi_t-\xi_t^{\pm})_{t\in\R}$ are bounded, one gets that
\be\label{tdelta}
\sup_{t\ge0}|\xi_{t+\delta}-\xi_t|<+\infty.
\ee\par
Finally, notice that $1\ge u(0,\cdot)\ge(2/3)\mathds{1}_{(-\infty,\xi^-_0)}$ a.e. in $\R$ by~\eqref{xi-t}, where $\mathds{1}_E$ denotes the characteris\-tic function of a measurable set $E$. The comparison principle~\eqref{comp} implies that, for every $t>0$,
$$1\ge u(t,\cdot)\ge v(t,\cdot)\hbox{ a.e. in }\R,$$
where $v$ denotes the solution of the Cauchy problem~\eqref{cauchy} with initial condition~$(2/3)\mathds{1}_{(-\infty,\xi^-_0)}$. It follows then from the spreading property~\eqref{spreading} that $\mathop{\essinf}_{(0,ct)}v(t,\cdot)\to1$ as $t\to+\infty$ for every $c>0$, whence
$$\mathop{\essinf}_{(0,ct)}u(t,\cdot)\to1\hbox{ as }t\to+\infty\hbox{ for every }c>0.$$
For any given $c>0$, property~\eqref{xi+t} then yields $\xi_t^+\ge ct$ for all $t>0$ large enough, whence $\liminf_{t\to+\infty}\xi_t^+/t\ge c$. Therefore,
$$\liminf_{t\to+\infty}\frac{\xi_t}{t}\ge c\ \hbox{ for all }c>0,$$
since $(\xi_t-\xi^+_t)_{t\in\R}$ is bounded. As a conclusion, $\xi_t/t\to+\infty$ as $t\to+\infty$, which contradicts~\eqref{tdelta}. The proof of Theorem~\ref{th1} is thereby complete.\hfill$\Box$


\SE{Stretching at large time for the Cauchy problem~\eqref{cauchy}}\label{sec3}

The second main result is concerned with flattening and stretching properties at large time for the solutions of the Cauchy problem~\eqref{cauchy}. We assume in this section that the equation~\eqref{eq} is homogeneous, in the sense that, for any $u_0\in L^{\infty}(\R;[0,1])$ and for any $h\in\R$, there holds
\be\label{homo}
u(t,\cdot+h)=v(t,\cdot)\hbox{ a.e. in }\R\hbox{ for every }t>0,
\ee
where $u$ and $v$ denote the solutions of the Cauchy problem~\eqref{cauchy} with initial conditions~$u_0$ and~$v_0:=u_0(\cdot+h)$. The property~\eqref{homo} is satisfied by the previous examples of dispersal operators~$\mathcal{D}u$: fractional Laplacian $\mathcal{D}u=-(-\Delta_x)^{\alpha}u$, the convolution operator $\mathcal{D}u=J*u-u$, the nonlinear fast diffusion $\mathcal{D}u=(u^\gamma)_{xx}$ and the nonlinear fractional diffusion $\mathcal{D}u=-(-\Delta_x)^{\alpha}(u^\gamma)$. If follows from assumptions~\eqref{comp} and~\eqref{homo} that if $u_0\in L^{\infty}(\R;[0,1])$ is nonincreasing, then
$$0\le u_0(\cdot+h)\le u_0\le1\hbox{ a.e. in }\R\hbox{ for every }h>0,$$
whence $0\le u(t,\cdot+h)\le u(t,\cdot)\le1$ a.e. in $\R$ for every $t>0$. In other words, $u(t,\cdot)$ is nonincreasing for every $t>0$. Under these conditions, for every $t\ge0$, denote
$$m(t)=\mathop{\essinf}_{\R}u(t,\cdot)=\mathop{\essinf}_{(x,+\infty)}u(t,\cdot)\ \hbox{ for every }x\in\R$$
and
$$M(t)=\mathop{\esssup}_{\R}u(t,\cdot)=\mathop{\esssup}_{(-\infty,x)}u(t,\cdot)\ \hbox{ for every }x\in\R,$$
and observe that $0\le m(t)\le M(t)\le 1$. Lastly, for every $\lambda\in(0,1)$, set
\be\label{xlambdat}
x_{\lambda}(t)=\inf\Big\{x\in\R;\ \mathop{\esssup}_{(x,+\infty)}u(t,\cdot)<\lambda\Big\}.
\ee
Notice that $x_{\lambda}(t)=+\infty$ if $\lambda\le m(t)$, while $x_{\lambda}(t)=-\infty$ if $\lambda>M(t)$. Moreover,
$$-\infty\le x_b(t)\le x_a(t)\le+\infty\hbox{ if }0<a\le b<1.$$

\begin{theo}\label{th2}
Assume that~\eqref{eq} satisfies~\eqref{comp},~\eqref{spreading} and~\eqref{homo} and let $u$ be the solution of the Cauchy problem~\eqref{cauchy} with a nonincreasing initial condition $u_0\in L^{\infty}(\R;[0,1])\backslash\{0\}$ such that $m(t)\to0$ and $M(t)\to1$ as $t\to+\infty$. Then, for every $0<a<b<1$,
\be\label{stretch}
\limsup_{t\to+\infty}\big(x_a(t)-x_b(t))=+\infty.
\ee
\end{theo}

Theorem~\ref{th2} means that the function $u(t,\cdot)$ becomes regularly as stretched as wanted as $t\to+\infty$, that is, for any given $0<\epsilon<1/2$, the ``reaction" zone where $u(t,\cdot)$ ranges between $\epsilon$ and $1-\epsilon$ cannot stay bounded at large time. However, the following reasonable conjecture remains open:

\begin{conj}\label{conj1}
Under the assumptions of Theorem~$\ref{th2}$, $x_a(t)-x_b(t)\to+\infty$ as $t\to+\infty$ for any $0<a<b<1$.
\end{conj}

The assumptions $m(t)\to0$ and $M(t)\to1$ are satisfied for instance if $\mathcal{D}u=-(-\Delta_x)^{\alpha}u$ with $0<\alpha<1$, $\mathcal{D}u=J*u-u$,  $\mathcal{D}u=(u^\gamma)_{xx}$ with $0<\gamma<1$, and $\mathcal{D}u=-(-\Delta_x)^{\alpha}(u^\gamma)$ with $0<\alpha<1$ and $\max(1-2\, \alpha,0)<\gamma$, and if $m(0)=0$ and $M(0)=1$ (remember that the $C^1([0,1])$ function~$f$ is such that $f(0)=f(1)=0$). These conditions are also satisfied for the same type of dispersal $\mathcal{D}u$ in the case $m(0)=0$, $0<M(0)\le1$ with $f>0$ on $(0,1)$.

Theorem~\ref{th2} is in spirit similar to a result obtained in~\cite{hr1} for the solutions of the Cauchy problem of the equation $u_t=u_{xx}+f(u)$: if $f$ satisfies~\eqref{f} and $f(u)/u$ is non-increasing over~$(0,1]$ and if the initial condition $u_0:\R\to(0,1]$ satisfies~$\liminf_{x\to-\infty}u_0(x)>0$, $u_0(+\infty)=0$, together with~$\lim_{x\to+\infty}\big(u_0(x)e^{\eta x}\big)=+\infty$ for every $\eta>0$ and $u'_0/u_0\in L^p(\R)\cap C^{2,\theta}(\R)$ for some~$p\in(1,+\infty)$ and $\theta\in(0,1)$, then $\|u_x(t,\cdot)\|_{L^{\infty}(\R)}\to0$ as $t\to+\infty$. This last conclusion is clearly stronger than~\eqref{stretch}. But, for our equation~\eqref{eq}, we make weaker regularity assumptions on $u_0$. In particular, $u_0$ may well be discontinuous and, since there is in general no spatial regulari\-zing effect for~\eqref{eq},~$u(t,\cdot)$ may not be differentiable or even continuous. But the property~\eqref{stretch} still gives some information about the profile of $u(t,\cdot)$ at large time. In particular, even if $u(t,\cdot)$ may be discontinuous at some points, the discontinuity jumps cannot stay large as time goes to $+\infty$.

Moreover, we can describe the profile of the solution $u$ around the position of each level sets. More precisely, for any level $\lambda\in(0,1)$, $u(t,\cdot+x_\lambda(t))$ becomes regularly as flat as wanted at large times and becomes as close as wanted to $\lambda$ locally on the left and on the right of the position $x_{\lambda}(t)$.

\begin{cor}\label{cor1}
Under the same assumptions as in Theorem~$\ref{th2}$, for every $0<\lambda<1$ and~$R>0,$ there holds
  \be\label{flat_moving}\left\{\baa{l}
\displaystyle\liminf_{t\to+\infty}\Big(\mathop{\esssup}_{x\in[-R,0]}\,|u(t,x+x_\lambda(t))-\lambda|\Big)=0,\vspace{3pt}\\
\displaystyle\liminf_{t\to+\infty}\Big(\mathop{\esssup}_{x\in[0,R]}\,|u(t,x+x_\lambda(t))-\lambda|\Big)=0.\eaa\right.
\ee
\end{cor}

\noindent{\bf{Proof of Theorem~\ref{th2}.}} Let $u$ be a solution of the Cauchy problem~\eqref{cauchy} with a nonincreasing initial condition $u_0\in L^{\infty}(\R;[0,1])\backslash\{0\}$. As already emphasized, $u(t,\cdot)$ is nonincreasing in $\R$ for every $t>0$. In particular, if $m(t)<\lambda<M(t)$, then $x_{\lambda}(t)$ defined in~\eqref{xlambdat} is a real number and
\be\label{xlambdatbis}
\mathop{\esssup}_{(x_{\lambda}(t),+\infty)}u(t,\cdot)\le\lambda\ \hbox{ and }\ \mathop{\essinf}_{(-\infty,x_{\lambda}(t))}u(t,\cdot)\ge\lambda.
\ee\par
Let now $a$ and $b$ be two given real numbers such that $0<a<b<1$ and assume by contradiction that $\limsup_{t\to+\infty}\big(x_a(t)-x_b(t)\big)<+\infty$, that is, there exist $T>0$ and $M>0$ such that
\be\label{xabt}
0\le x_a(t)-x_b(t)\le M\ \hbox{ for all }t\ge T.
\ee
Since $m(t)\to0$ and $M(t)\to1$ as $t\to+\infty$, one can assume without loss of generality that~$m(t)<a<b<M(t)$ for all $t\ge T$, whence $x_a(t)$ and $x_b(t)$ are real numbers for $t\ge T$. Notice that
\be\label{ab2}
x_b(t)\le x_{(a+b)/2}(t)\le x_a(t)\ \hbox{ for all }t\ge T.
\ee
From the uniform continuity of $u$ with respect to $t$ in $[T,+\infty)\times\R$, there is $\delta>0$ such that, for every $(t,s,x)\in[T,+\infty)\times[T,+\infty)\times\R$,
$$\big(|t-s|\le\delta\big)\ \Longrightarrow\ \Big(|u(t,x)-u(s,x)|\le\frac{b-a}{4}\Big).$$
Now, for every $t\ge T$, one has $u(t,\cdot)\ge b$ a.e. in $(-\infty,x_b(t))$, whence
$$u(t+\delta,\cdot)\ge b-\frac{b-a}{4}=\frac{3b+a}{4}\hbox{ a.e. in }(-\infty,x_b(t)),$$
whereas $u(t,\cdot)\le a$ a.e. in $(x_a(t),+\infty)$ and
$$u(t+\delta,\cdot)\le a+\frac{b-a}{4}=\frac{3a+b}{4}\hbox{ a.e. in }(x_a(t),+\infty).$$
Since $(3a+b)/4<(a+b)/2<(3b+a)/4$, one infers from the definition of $x_{(a+b)/2}(t+\delta)$ that
$$x_b(t)\le x_{(a+b)/2}(t+\delta)\le x_a(t),$$
whence
\be\label{xab2bis}
\big|x_{(a+b)/2}(t+\delta)-x_{(a+b)/2}(t)\big|\le M\ \hbox{ for all }t\ge T.
\ee
from~\eqref{xabt} and~\eqref{ab2}.\par
On the other hand, since $u(T,\cdot)$ is nonincreasing and ranges in $[0,1]$, it follows again from the definition of $x_{(a+b)/2}(T)$ that
$$1\ge u(T,\cdot)\ge\frac{a+b}{2}\mathds{1}_{(-\infty,x_{(a+b)/2}(T))}=:v_0\ \hbox{ a.e. in }\R.$$
Therefore, the comparison principle~\eqref{comp} implies that, for every $t>0$, $1\ge u(t+T,\cdot)\ge v(t,\cdot)$ a.e. in $\R$, where $v$ denotes the solution of the Cauchy problem~\eqref{cauchy} with initial condition~$v_0\in L^{\infty}(\R;[0,1])\backslash\{0\}$. The assumption~\eqref{spreading} applied to $v$ therefore yields $\essinf_{(0,ct)}v(t,\cdot)\to1$ as $t\to+\infty$ for every $c>0$, whence
$$\mathop{\essinf}_{(0,ct)}u(t+T,\cdot)\to1\hbox{ as }t\to+\infty,\hbox{ for every }c>0.$$
Finally, for any arbitrary positive real number $c$, there holds
$$\mathop{\essinf}_{(0,cn\delta)}u(T+n\delta,\cdot)>\frac{a+b}{2}$$
for $n\in\N$ large enough, whence $x_{(a+b)/2}(T+n\delta)\ge cn\delta$ for $n$ large enough, since $u(T+n\delta,\cdot)$ is nonincreasing. As a consequence,
$$\liminf_{n\to+\infty}\frac{x_{(a+b)/2}(T+n\delta)}{n\delta}\ge c$$
and, since $c>0$ is arbitrary, it follows that $x_{(a+b)/2}(T+n\delta)/(n\delta)\to+\infty$ as $n\to+\infty$. This contradicts~\eqref{xab2bis} and the proof of Theorem~\ref{th2} is thereby complete.\hfill$\Box$

\begin{rem}{\rm Consider here the case of equation~\eqref{fraclaplace} with the fractional Laplacian $\mathcal{D}u=-(-\Delta_x)^{\alpha}u$ with $0<\alpha<1$. For this equation, there is a regularizing effect in the spatial variable and, under the assumptions of Theorem~\ref{th2}, it follows that the solution $u$ is of class~$C^1$ in $[\delta,+\infty)\times\R$ for every $\delta>0$, with bounded derivatives with respect to $t$ and~$x$. For any given $\lambda\in(0,1)$, there is a real number $T_{\lambda}>0$ such that, for every $t\ge T_{\lambda}$, there is a unique~$x_{\lambda}(t)\in\R$ such that $u(t,x_{\lambda}(t))=\lambda$. We now claim that
$$\liminf_{t\to+\infty}|u_x(t,x_{\lambda}(t))|=0.$$
Indeed, otherwise, there are $\epsilon>0$ and $T'\ge T_{\lambda}$ such that $u_x(t,x_{\lambda}(t))\le-\epsilon$ for all $t\ge T'$. It follows then from the implicit function theorem that the function~$x_{\lambda}$ is of class $C^1([T',+\infty))$ and that
$$x_{\lambda}'(t)=-\frac{u_t(t,x_{\lambda}(t))}{u_x(t,x_{\lambda}(t))}$$
for all $t\ge T'$. Therefore, there is $M\ge0$ such that $|x_{\lambda}'(t)|\le M$ for all $t\ge T'$ and one can get a contradiction as in the end of the proof of Theorem~\ref{th2}.}
\end{rem}

Again for~\eqref{fraclaplace} with the fractional Laplacian $\mathcal{D}u=-(-\Delta_x)^{\alpha}u$, and with a Fisher-KPP function $f$, Roquejoffre and Tarfulea~\cite{RoqTar15} recently showed that the partial derivatives of $u$ with respect to $x$ are all bounded by~$u$ times an exponentially decaying function, that is, for every $k\in\N$, there are positive constants $C_k$ and $\delta_k$ such that $\|\partial^k_x u(t,\cdot)\|_\infty\leq C_k u(t,x) e^{-\delta_k t}$ for large $t$. This result has been proved under the additional condition that~$u_0$ is exponentially decaying as $|x|\to+\infty$, and it actually holds in any spatial dimension. It implies in particular that $\|u_x(t,\cdot)\|_{L^{\infty}(\R)}\to0$ as $t\to+\infty$. Moreover, it is proved in~\cite{RoqTar15} that the large time dynamics of the KPP problem~\eqref{fraclaplace} is in some sense the same as that of the ODE $\theta'(t)=f(\theta).$

More generally speaking, under the general conditions~\eqref{comp},~\eqref{spreading} and~\eqref{homo}, we conjecture that fast propagation should lead to a flattening of the solution in the sense that, at least if $f>0$ on $(0,1)$, then the only $x$-monotone solutions $u:\R\times\R\to(0,1)$ of~\eqref{fraclaplace} are actually $x$-independent, that is, they can be written as $u(t,x)=\theta(t)$, where $\theta'(t)=f(\theta(t))$ for all $t\in\R$ with $\theta(-\infty)=0$ and $\theta(+\infty)=1$.

An intermediate result that we prove here under the sole conditions~\eqref{comp},~\eqref{spreading} and~\eqref{homo} is the flattening result stated in Corollary~\ref{cor1}. We point out that this result holds for all equations~\eqref{fraclaplace},~\eqref{convolution},~\eqref{eq:fastdif} and~\eqref{eq:fastdif_frac} under the same conditions as in Section~\ref{intro}.\hfill\break

\noindent{\bf{Proof of Corollary~\ref{cor1}.}} Let $u$ be a solution of the Cauchy problem~\eqref{cauchy} with a nonincreasing initial condition $u_0\in L^{\infty}(\R;[0,1])\backslash\{0\}$. As already emphasized, $u(t,\cdot)$ is nonincreasing in $\R$ for every $t>0$. In particular, if $m(t)<\lambda<M(t)$, then $x_{\lambda}(t)$ defined in~\eqref{xlambdat} is a real number and~\eqref{xlambdatbis} holds.\par
Let then $\lambda$ be any real number such that $0<\lambda<1,$ $R$ be any positive real number and $\epsilon>0$ be any positive real number such that $0<\lambda-\epsilon<\lambda<\lambda+\epsilon<1$. From Theorem~\ref{th2}, there exist two increasing sequences $(t_n)_{n\in\N}$ and $(t'_n)_{n\in\N}$ of positive real numbers such that
$$x_{\lambda}(t_n)-x_{\lambda+\epsilon}(t_n) \to +\infty\ \hbox{ and }\ x_{\lambda-\epsilon}(t'_n)-x_{\lambda}(t'_n) \to +\infty\ \hbox{ as } n\to+\infty.$$
As a consequence, there exists $N\in\N$ such that, for all $n\ge N$,
$$\left\{\baa{ll}
x_{\lambda+\epsilon}(t_n)\leq  x+x_\lambda(t_n) & \hbox{for all }x\in[-R,0],\vspace{3pt}\\
x+x_\lambda(t'_n)\leq x_{\lambda-\epsilon}(t'_n) & \hbox{for all }x\in[0,R].\eaa\right.$$
Thus, by~\eqref{xlambdatbis} applied with $\lambda-\epsilon$, $\lambda$ and $\lambda+\epsilon$, one infers that, for all $n\ge N$,
$$\left\{\baa{ll}
\lambda\le u(t_n,x+x_{\lambda}(t_n))\le\lambda+\epsilon & \hbox{for a.e. }x\in[-R,0],\vspace{3pt}\\
\lambda-\epsilon\le u(t'_n,x+x_{\lambda}(t'_n))\le\lambda & \hbox{for a.e. }x\in[0,R],\eaa\right.$$
whence
$$\mathop{\esssup}_{x\in[-R,0]}|u(t_n,x+x_\lambda(t_n))-\lambda|\leq\epsilon\ \hbox{ and }\ \mathop{\esssup}_{x\in[0,R]}|u(t'_n,x+x_\lambda(t'_n))-\lambda|\leq\epsilon$$
for all $n\geq N$. Since $\epsilon>0$ can be arbitrarily small and $\lim_{n\to+\infty}t_n=\lim_{n\to+\infty}t'_n=+\infty$, the proof of Corollary~\ref{cor1} is thereby complete.\hfill$\Box$

\begin{rem}{\rm It easily follows from the proof of Corollary~\ref{cor1} that, if Conjecture~\ref{conj1} is true, then for any level $\lambda\in(0,1),$ the solution $u$ of the Cauchy problem~\eqref{cauchy} converges to $\lambda$ locally around each position of the level set as $t\to+\infty$, that is,
$$\mathop{\esssup}_{x\in[-R,R]}|u(t,x+x_\lambda(t))-\lambda|\to0\ \hbox{ as }t\to+\infty$$
for every $R>0$. This last convergence result would then hold automatically for the whole family $t\to+\infty$ and on both the left and right of $x_{\lambda}(t)$ simultaneously.}
\end{rem}

\begin{figure}
\centering
\subfigure[]{\includegraphics[width=0.7\textwidth]{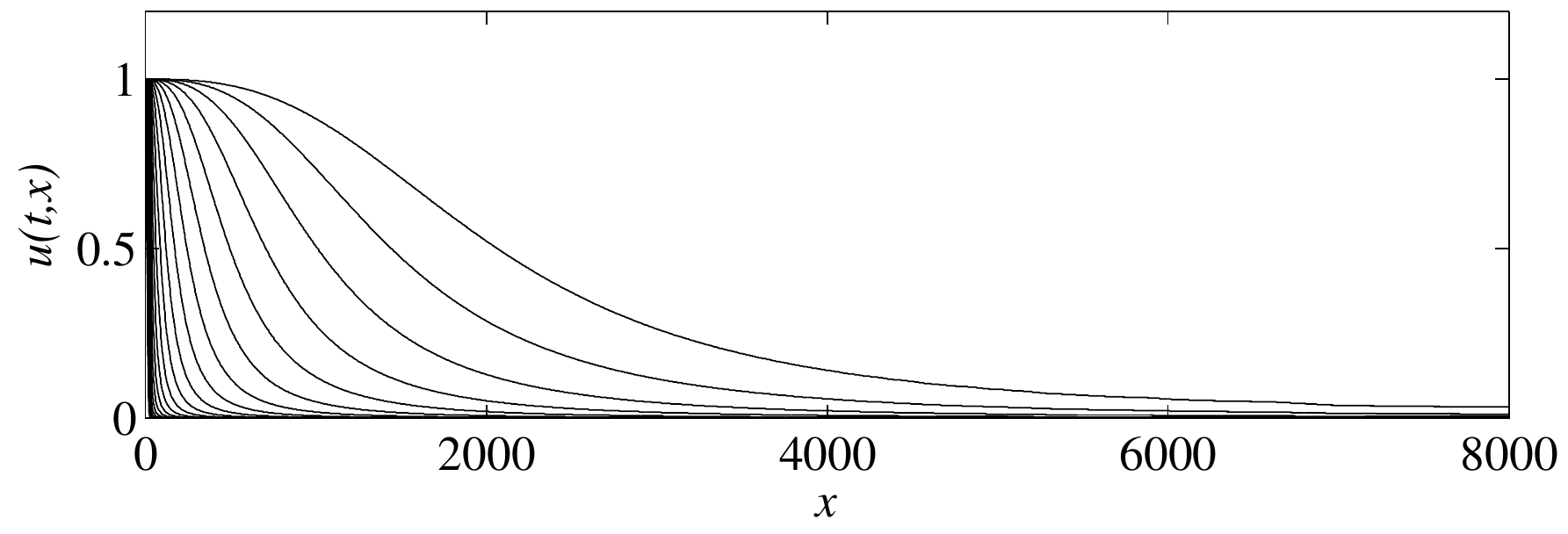}}
\subfigure[]{\includegraphics[width=0.7\textwidth]{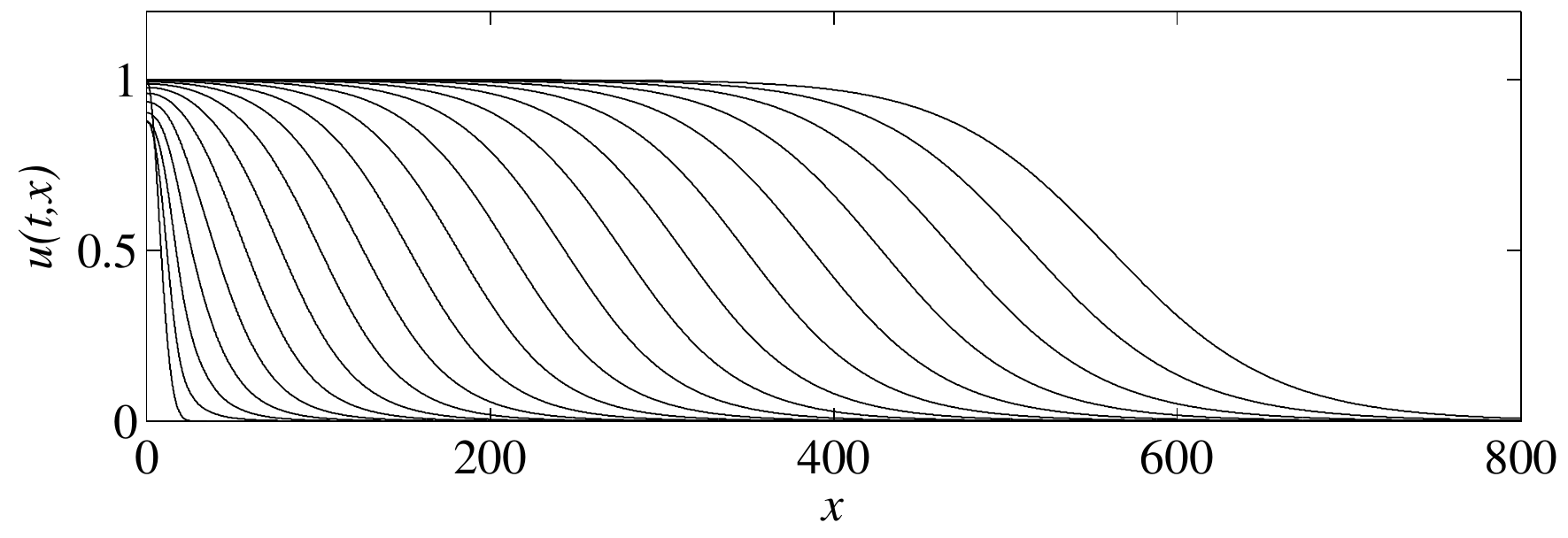}}
\subfigure[]{\includegraphics[width=0.7\textwidth]{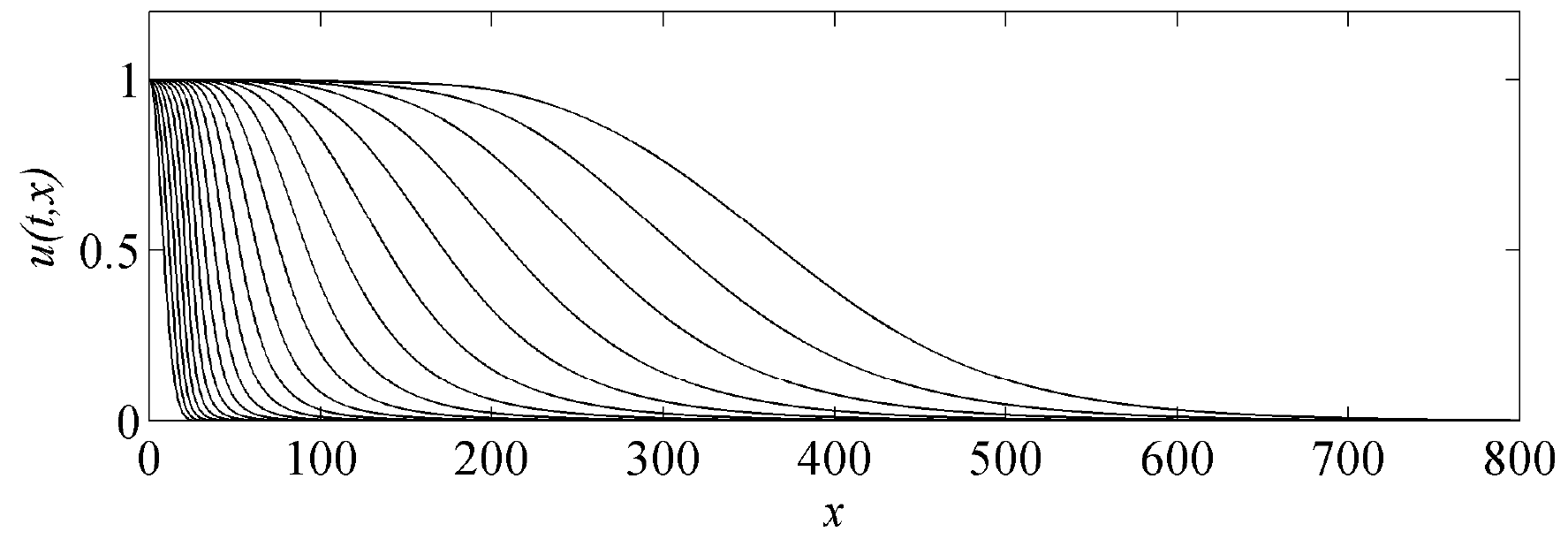}}
\subfigure[]{\includegraphics[width=0.7\textwidth]{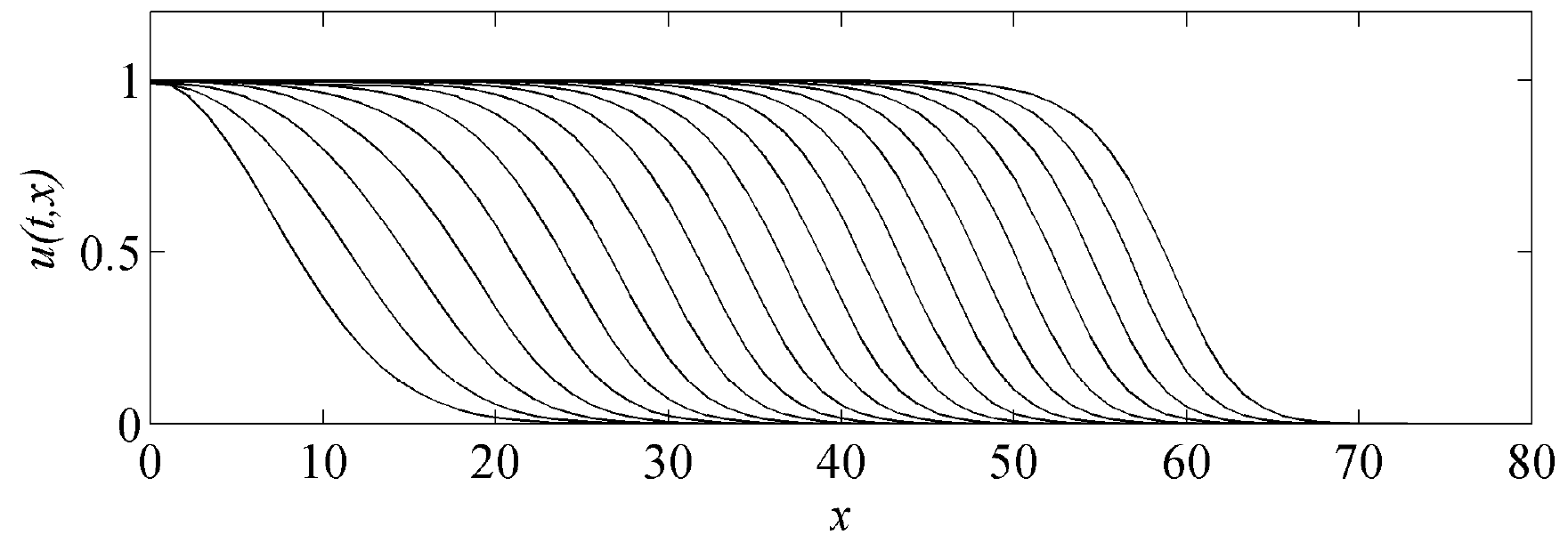}}
\caption{Propagation to the right of the solution of the Cauchy problem \eqref{cauchy}, at successive times $t=0,1,\ldots,20,$ with  (a) $\mathcal D u = -(-\Delta_x)^\alpha u$ with $\alpha=0.9$; (b)  $\mathcal D u = J*u-u$ with $J(x)=\exp(-\sqrt{|x|})/4$;  (c) $\mathcal D u=(u^\gamma)_{xx}$ with $\gamma=1/2$; and (d) $\mathcal D u = u_{xx}$. In all cases, the initial condition was $u_0(x)=\exp(-x^2/100)$ and the function $f$ was of the KPP type $f(u)=u\, (1-u)$.}
 \label{fig:1}
\end{figure}

\begin{figure}
\centering
\includegraphics[width=0.4\textwidth]{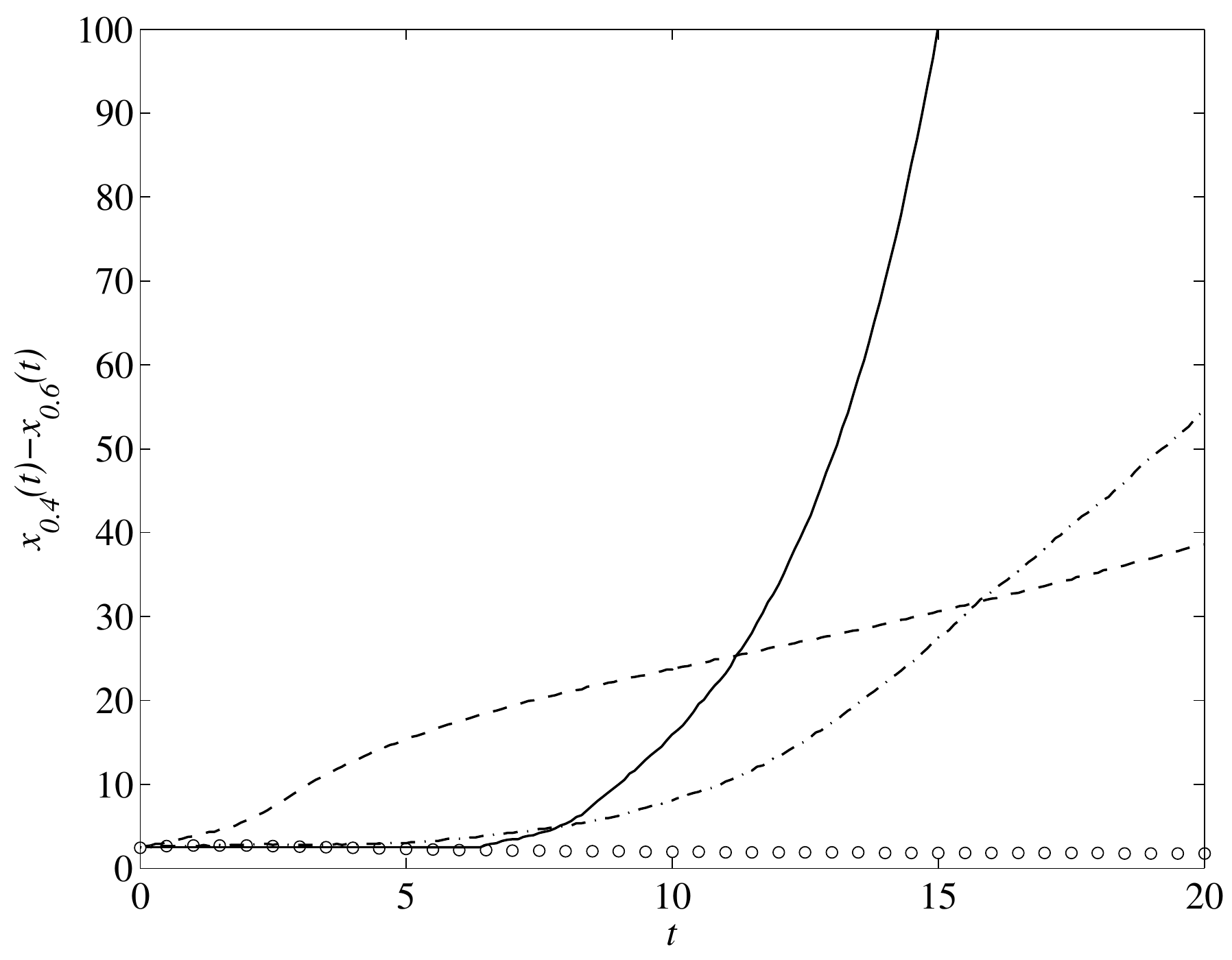}
\caption{Distance $x_{0.4}(t)-x_{0.6}(t)$ between two level sets of the solution of the Cauchy problem \eqref{cauchy}, for $t\in(0,20)$, with: (solid line) $\mathcal D u = -(-\Delta_x)^\alpha u$ with $\alpha=0.9$; (dashed line)  $\mathcal D u = J*u-u$ with $J(x)=\exp(-\sqrt{|x|})/4$;  (dash-dot line) $\mathcal D u=(u^\gamma)_{xx}$ with $\gamma=1/2$; and (circles) $\mathcal D u = u_{xx}$. The assumptions on $f$ and $u_0$ are the same as in Fig.~\ref{fig:1}.}
 \label{fig:2}
\end{figure}

Lastly, in order to illustrate the theoretical results on the flattening and stretching properties of the solutions of the Cauchy problem \eqref{cauchy}, we have performed some numerical simulations that are depicted in figures~\ref{fig:1} and~\ref{fig:2}. They have been done with the KPP nonlinearity $f(u)=u\, (1-u),$ for the three examples $\mathcal D u = -(-\Delta_x)^\alpha u$ with $\alpha=0.9$,  $\mathcal D u = J*u-u$ with $J(x)=\exp(-\sqrt{|x|})/4$, and  $\mathcal D u=(u^\gamma)_{xx}$ with $\gamma=1/2$ (Fig.\ref{fig:1},~a,~b and~c, respectively). These properties, and the acceleration of the solutions in each of these three cases are to be compared with the well-known convergence to a standard travelling front which can be observed with the standard diffusion operator $\mathcal D u = u_{xx}$ (Fig.~\ref{fig:1}, d). To compute these solutions, we used (a and b): the Strang splitting, which consists in splitting the equation \eqref{cauchy} into two simpler evolution problems \cite{Cou14}: $v_t =-(-\Delta_x)^\alpha v$ and  $v_t =J*v-v$ which are treated with fast Fourier transform techniques and $w_t =w\, (1-w)$ which solution can be computed explicitly; (c and d): the time-dependent finite element solver of Comsol Multiphysics$^{\copyright}$.

\end{document}